\documentclass[12pt]{article}
\usepackage{amsmath,amssymb,amsfonts,amsthm}
\newenvironment{myabstract}{\par\noindent
{\bf Abstract . } \small }
{\par\vskip8pt minus3pt\rm}
\newcounter{item}[section]
\newcounter{kirshr}
\newcounter{kirsha}
\newcounter{kirshb}
\newenvironment{enumroman}{\setcounter{kirshr}{1}
\begin{list}{(\roman{kirshr})}{\usecounter{kirshr}} }{\end{list}}
\newenvironment{enumarab}{\setcounter{kirshb}{1}

\begin{list}{(\arabic{kirshb})}{\usecounter{kirshb}} }{\end{list}}
\newenvironment{athm}[1]{\vskip3mm\par\noindent
{\bf #1 }. \slshape }
{\upshape\par\vskip10pt minus3pt}
\newenvironment{demo}[1]{\noindent{\bf #1.}\upshape\mdseries}
{\nopagebreak{\hfill\rule{2mm}{2mm}\nopagebreak}\par\normalfont}
\theoremstyle{definition}

\def\Sg{{\mathfrak{Sg}}}

\def\CA{{\bf CA}}
\def\RCA{{\bf RCA}}

\def\Rd{{\ Rd}}

\def\(R)RA{{\bf (R)RA}}

\def\R{\mathbb{R}}

\def\C{\mathbb{C}}
\def\A{{\mathfrak{A}}}
\def\B{{\mathfrak{B}}}
\def\C{{\mathfrak{C}}}

\def\R{{\mathfrak{R}}}
\def\P{{\mathfrak{P}}}

\def\Rd{{\mathfrak{Rd}}}

\def\A{{\mathfrak{A}}}
\def\B{{\mathfrak{B}}}
\def\C{{\mathfrak{C}}}

\def\Rd{{\mathfrak{Rd}}}

\def\RQA{{\bf RQA}}

\def\RQEA{{\bf RQEA}}

\title{The class of representable quasipolyadic equality algebras of dimension $\omega$ 
is not finitely axiomatizable over the class of representable quasipolyadic algebras of the same dimension }
\begin{document}
\maketitle

\begin{myabstract} Modifying Andr\'eka's methods of splitting, we show that the class $\RQEA_{\omega}$ 
of representable quasipolyadic equality algebas of dimension $\omega$ is not axiomatized by a set of universal formulas 
containing only finitely many variables over the class $\RQA_{\omega}$ of representable quasipolyadic algebras of dimension
$\omega$. 
This answers a question posed by  Andr\'eka and generalizes a result of hers proved for finite dimensions $>2$.

\end{myabstract}

\section{Introduction}

Stone's representation theorem for Boolean algebras can be formulated in two, essentially equivalent ways. 
Every Boolean algebra is isomorphic to a field of sets,
or the class of Boolean set algebras can be axiomatized by a finite set of equations. As is well known, Boolean algebras 
constitute the algebraic counterpart of propositional logic.
Stone's representation theorem, on the other hand, is the algebraic equivalent of the completeness theorem for propositional logic.

However, when we step inside the realm of first order logic, things tend to become more complicated. 
Not every abstract cylindric algebra is representable as a field of sets, where
the extra Boolean operations of cylindrifiers and diagonal elements are faithfully represented by projections and equality.
Disappointingly, the class of representable algebras fail to be axiomatized  by any reasonable finite schema and its resistance to such 
axiomatizations is inevitable. This is basically a reflection of the essential incompleteness of natural (more basic) infinitary 
extensions of first order logic. In such extensions, unlike first order logic, validity cannot be captured by a finite schema. 

Such extentions are obtained by dropping the condition of local finiteness 
(reflecting the simple fact that first order formulas contain only finitely many variables)
in algebras considered, 
allowing formulas of infinite length. This is necessary if we want to deal with the so-called algebriasable extensions of first order logic; 
extensions that are akin to universal algebraic investigations.

The condition of local finiteness, cannot be expressed in first order logic, 
and this is not warranted if we want to deal, like in the case of Boolean algebras, only with equations, 
or at worst quasi-equations.  Then we are faced with the following problem. Find a simple (hopefully finitary) axiomatization 
of classes of representable algebras abounding in algebraic logic, using only equations or quasi equations, 
which also means that we want to stay in the realm of quasivarieties.

There are  two conflicting but complementary facets 
of such a problem, referred to in the literature, as the representation problem.
One is to delve deeply in investigating the complexity of potential axiomatizations for existing varieties 
of representable algebras, which we do in this paper, the other is to try to sidestep 
such wild unruly complex axiomatizations, often referred to as {\it taming methods}. 
Those taming methods can either involve passing to (better behaved) expansions of the algebras considered,
or else change the very  notion of representatiblity involved, as long as it remains concrete enough.
The borderlines are difficult to draw, we do might not know what is {\it not} concrete enough, but we can
 judge that a given representability notion is satisfactory, once we have one. 
(This is analogous to  undecidability issues, with the main difference that we do know what we mean 
by {\it not decidable}. We do not have an analogue of a 'recursive representability notion').

The representation problem is fairly well understood for cylindric algebras, however, 
for some of its relatives it is not as deeply explored.

Quoting Henkin Monk and Tarski in \cite{HMT2} p.266-267:
``Quasipolyadic algebras: These are like polyadic algebras, except that ${\sf s}_{\tau}$ is allowed only for finite transformations, 
and ${\sf c}_{(\Gamma)}$ only for finite
$\Gamma$. Their theory has not been much developed, but they form an interesting stage between cylindric and polyadic algebras"
The results in this paper are  a contribution to the representation theory of quasipolyadic (equality) algebras, that has recently been investigated in the 
literature, coming to the forefront of algebraic logic again with a plathora of new results obtained by Andr\'eka, N\'emeti, Ferenczi, 
Hodkinson, S\'agi, the present author and others, cf. \cite{Fer}, \cite{F1}, \cite{F2}, \cite{F3}, \cite{ES}, \cite{amal}, \cite{sf}, \cite{h} 
\cite{s}, \cite{Sagi}, \cite{ns}, \cite{ANS}, \cite{c}. 

Let $U$ be a set and $\alpha$ be a ordinal. Then $\B(^{\alpha}U)$ is the boolean set algebra 
with unit $^{\alpha}U$.
Let $\tau:\alpha\to \alpha,$ $i,j<\alpha$
and $X\subseteq {}^{\alpha}U.$ 
Then 
$${\sf s}_{\tau}^UX=\{s\in {}^{\alpha}U: s\circ \tau\in X\},$$
$${\sf c}_i^UX=\{t\in {}^{\alpha}U: \exists s\in X \text { and } t(j)=s(j) 
\text { for all }j\neq i\},$$
and $${\sf d}_{ij}^U=\{s\in {}^{\alpha}U: s_i=s_j\}.$$
Superscripts are omitted if no confusion is likely to ensue.
${\bf SP}$ stands for the operation of forming subdirect products.
$$\RCA_{\alpha}=
{\bf SP}\{\langle \B(^{\alpha}U), 
{\sf c}_i^U, {\sf d}_{ij}^U \rangle_{i,j<\alpha}:  U \text { is a set }\}$$
It is known \cite{Andreka}  that for $\alpha\geq 3$, the class $\RCA_{\alpha}$ of $\alpha$ dimensional representable cylindric algebras 
is a variety that is not axiomatizable by any set $\Sigma$ 
of quantifier free formulas containing only finitely
many free variables. Let $$\RQA_{\alpha}=
{\bf SP}\{\langle \B(^{\alpha}U), 
{\sf c}_i^U,  {\sf s}_{[i|j]}^U, {\sf s}_{[i,j]}^U\rangle_{i,j<\alpha}:  U \text { is a set }\}.$$
Let $$\RQEA_{\alpha}=
{\bf SP}\{\langle \B(^{\alpha}U), 
{\sf c}_i^U, {\sf d}_{ij}^U, {\sf s}_{[i,j]}^U\rangle_{i,j<\alpha}:  U \text { is a set }\}.$$
The map $[i,j]$ is the tansposition that interchanges $i$ and $j$ that is $[i,j](i)=j$ and $[i,j](j)=i$ and
$[i,j]x=x$ for $x\notin \{i,j\}$. On the other hand, the map $[i|j]$ is the replacement that sends $i$ to $j$ and is the identity otherwise.  
In $\RCA_{\alpha}$ and $\RQEA_{\alpha}$, ${\sf s}_{[i|j]}$ is term definable by ${\sf c}_i(x\cdot {\sf d}_{ij})$ for $i\neq j$; 
according to a widespread custom, we denote ${\sf s}_{[i|j]}$ by ${\sf s}_j^i.$ The term ${\sf s}_i^ix$ is just $x$.

It is also known that $\RQEA_{\alpha}$ is a variety \cite{ST}, that is it is closed under homomorphic images.
It seems to be a  known result (at least implicitly) that for $\alpha$ an infinite ordinal, $\RQEA_{\alpha}$ cannot be axiomatized 
by a set of universal formulas that contains only finitely many variables.
However, no complete proof of the latter result exists in the literature (a sketch of proof is given in \cite{ST} but it has serious gaps, see below.) 
The proof for finite $\alpha$ is in \cite{Andreka}. 
We strengthen this result by showing that $\RQEA_{\omega}$ is not even finitely axiomatizable over $\RQA_{\omega}$. That is, 
we show that the diagonal elements cannot be captured by any set of universal formulas containing finitely many 
variables. This generalizes a result of Andr\'eka proved for the finite dimensional case in \cite{Andreka}.

Sometimes the theory of quasipolyadic equality algebras 
coincide with that of cylindric algebras like in the cases of finite dimensional algebras of positive characteristic, 
locally finite and dimension complemented 
algebras \cite{HMT2} of infinite dimensions; these are all representable algebras (isomorphic to subdirect products of set algebras),  
in which substitutions are term definable, 
but in general substitutions are {\it not} term definable in (even representable) 
cylindric algebras. Even more, 
it is also known \cite{ANS} that there is an infinite dimensional quasipolyadic equality algebra 
that is not representable with a representable cylindric reduct, that is the representability of the  cylindric reduct of a quasipolyadic equality algebra
does not force representability of substitutions concretely. This result can be strengthened to show that the ${\sf s}_{[i,j]}$'s are not even 
finitely axiomatizable over 
$\RCA_{\omega}$ as proved in \cite{t}.
This means that, on the intuitive level, the operations of substitutions 
add a lot to cylindric algebras, 
and suggests that generalizing results from cylindric algebras to quasipolyadic equality algebras, is far from a direct procedure, 
as indeed the proof of our main result herein, compared to Andr\'eka's, illustrates.

In our proofs, we modify  the techniques of splitting due to  Andr\'eka, which she applies to cylindric algebras, to adapt to 
the quasipolyadic equality case.
For finite dimensions this can be done relatively easily because the splitting is implemented relative to a finite set of substitutions 
and this involves a combinatorial trick depending on counting the number of substitutions. 
This technique no longer holds in the presence of infinitely many substitutions because  we simply cannot count them. 
So the main idea in this paper is that we do the splitting relative to reducts containing only finitely many substitutions and 
our desired algebra, witnessing the complexity of 
axiomatizations, will be a limit 
of such reducts. The idea might be simple, but the  details of implementing the idea turn out to be somewhat intricate; one 
has to make sure that certain (basically) combinatorial 
properties holding for all finite reducts pass to their limit.

The subtle splitting technique invented by Andr\'eka can be summarized as follows. 
In the presence of only finitely many substitutions, we take a fairly simple representable algebra generated by an atom, and we break 
up or {\it split} the atom into enough (finitely many) atoms,
forming a larger algebra, that is in fact non-representable; in fact, its cylindric reduct will not be representable, due to the  incompatibility 
between the number of atoms, and the number of elements in the domain of a representation. However, the ''small" subalgebras 
of such an algebra will be 
representable.
\footnote{This has affinity to Monk's construction of what is known as 'Monk's algebras' witnessing non finite axiomatizability for 
the class of representable cylindric algebras.
The key idea of the construction of a Monk's algebra is not so hard. 
Such algebras are finite, hence atomic, more precisely their Boolean reducts are atomic.
The atoms are given colours, and 
cylindrifications and diagonals are defined by stating that monochromatic triangles 
are inconsistent. If a Monk's algebra has many more atoms than colours, 
it follows from Ramsey's Theorem that any representation
of the algebra must contain a monochromatic traingle, so the algebra
is not representable.}

In our construction, though, we perform infinitely many finite splittings 
(not just one which is done in \cite{ST}), increasing in number but always finite, 
constructing infinitely many algebras, whose similarity types contain only finitely many 
substitutions. This is the main novelty occuring here, {\it  a modification of Andr\'eka's method of splitting to adapt to 
the quasipolyadic equality case.} Such constructed non-representable algebras, 
form a chain, and our desired algebra will be their 
directed union. The easy thing to do is to show that ``small" subalgebras of every non-representable 
algebra in the chain is representable; the hard thing to do is to show that ``small" subalgebras of the non-representable 
limit remain 
representable. (The error in Sain's Thompson paper is claiming that the small subalgebras 
of the non-representable algebra, obtained by performing only one splitting into infinitely many atoms, are representable; this is not necessarily 
true).

The cylindric reduct of the algebras forming the chain is of $\CA_{\omega}$ type; in particular, it contains infinitely many 
cylindrifications and diagonal elements.
The  combinatorial argument of counting depends essentially on the presence of infinitely many diagonal elements.
Indeed, it can be shown that the splitting 
technique adopted to prove complexity results concerning axiomatizations of $\RQEA_{\omega}$ 
simply does not work in the absence of diagonals. This can be easily destilled from our proof since our constructed 
non-representable  quasipolyadic equality algebras, in fact have a representable 
quasipolyadic reduct. 
An open problem here, that can be traced back to to Sain's and Thompson's paper \cite{ST}, 
is whether $\RQA_{\omega}$ {\it can be} 
axiomatized by a necessarily infinite) set of formulas using only finitely many 
variables. This seems to be a hard problem, and the author tends to believe 
that there are axiomatizations that contain only finitely many variables, but further research is needed in this area.

On the other hand, the algebra constructed by this method of splitting is `almost representable', in the sense that if we 
enlarge the potential domain of a representation, then various reducts of the algebra, obtained by discarding some of the operations 
(for example diagonal elements or infinitely many cylindrifications), turn out 
representable; and this gives {\it relative} non-finitizability results, witness Theorem 3 below and  \cite{t}.
Here we are encountered by a situation where we cannot have our cake and eat. If we want a quasipolyadic equality 
algebras that is only barely representable,
then we cannot obtain non-representability of some of its strict reducts like its quasipolyadic reduct.

Throughout, we will be tacitly assuming that quasipolyadic (equality) algebras are not only 
term-definitionally with finitary polyadic (equality) algebras as proved in \cite{ST} p.546, but that they are actually
the same. This means that in certain places we consider only substitutions 
corresponding to transpositions rather then all substitutions corresponding to finite transformations which is perfectly 
legitimate. Also we understand representability of reducts of quasipolyadic equality algebras, when we discard some of the substitution operations, 
in the obvious sense. 

\section{Main result and its proof}

If $\A$ has a cylindric reduct, then $\Rd_{ca}\A\in \RCA_{\alpha}$ denotes this reduct. 
Our next theorem corrects the error mentioned above in Sain's Thompson's seminal paper \cite{ST},
generalizes Theorem 6 in \cite{Andreka} p. 193 to infinitely many dimensions, 
and answers a question by Andre\'ka in op cit also on p. 193. 

\begin{athm}{Theorem 1} The variety $\RQEA_{\omega}$ cannot be axiomatized with with a set $\Sigma$ of quantifier free
formulas containing finitely many variables. In fact, for any $k<\omega$, and any set of quantifier free formulas 
$\Sigma$ axiomatizing $\RQEA_{\omega}$, $\Sigma$ contains a formula with more than $k$ variables in which some diagonal 
element occurs.
\end{athm}

\begin{demo}{Proof}  The proof consists of two parts.  In the first part we construct algebras $\A_{k,n}$ with certain properties, for each 
$n,k\in \omega\sim \{0\}$. In the second part we form a limit of such algebras as $n$ tends to infinity, obtaining an algebra $\A_k$ that is not representable, though its 
$k$-generated subalgebras are representable. This algebra wil finish the proof.

{\bf Part I}

Let $k, n\in \omega\sim \{0\}$.
Let $G_n$ be the symmetric group on $n$.
$G_n$ is generated by the set of all transpositions $\{[i,j]: i,j\in n\}$ and for $n\leq m$, we can consider $G_n\subseteq G_m$.
We shall construct an algebra $\A_{k,n}=(A_{k,n}, +, \cdot ,-, {\sf c}_i, {\sf s}_{\tau}, {\sf d}_{ij})_{i,j\in \omega, \tau\in G_n}$ 
with the following properties.
\begin{enumroman}
\item $\Rd_{ca}\A_{k,n}\notin \RCA_{\omega}$. 
\item Every $k$-generated subalgebra of $\A_{k,n}$ is representable.
\item 
There is a one to one mapping $h:\A_{k,n}\to (\B(^{\omega}W), {\sf c}_i, {\sf s}_{\tau}, {\sf d}_{ij})_{i,j<\omega, \tau\in G_n}$ 
such that $h$ is a homomorphism with respect to all operations of $\A_{k,n}$ except for the 
diagonal elements. 
\end{enumroman}

Here $k$-generated means generated by $k$ elements.
The proof for finite reducts uses arguments very similar to the proof of Andr\'eka of Theorem 6 in \cite{Andreka}, and has affinity  with the proof 
of theorem 3.1 in \cite{c}. 
However, there are two major differences.
Our cylindric reducts are infinite dimensional, and our proof is more direct and, in fact, far easier to grasp. 
The proof of the above cited theorem of Andr\'eka's goes through the route of certain finite expansions by so-called 
permutation invariant unary operations that are also modalities (distributive over the boolean join), and these are more general 
than substitutions. Substitutions are more concrete, and therefore our proof is less abstract.
\begin{enumarab}
\item Let $m\geq 2^{k.n!+1}$, $m<\omega$ 
and let $\langle U_i:i<\omega\rangle$ be a system
of disjoint sets such that $|U_i|=m$ for $i\geq 0$ and $U_0=\{0,\ldots m-1\}$.
Let $$U=\bigcup\{U_i: i\in \omega\},$$
let $$R=\prod _{i<\omega}U_i=\{s\in {}^{\omega}U: s_i\in U_i\},$$
and let $\A'$ be the subalgebra of 
$\langle \B({}^{\omega}U), {\sf c}_i, {\sf d}_{ij}, {\sf s}_{\tau}\rangle_{\tau\in G_n}$
generated by $R$.
Then ${\sf s}_{\tau}R$ is an atom of $\A'$ for any $\tau\in G_n$. Indeed for any two sequences $s,z\in R$ there is a 
permutation $\sigma:U\to U$ of $U$ taking $s$ to $z$ and fixing $R$, i.e
$\sigma\circ s=z$ and $R=\{\sigma\circ p:p\in R\}$. 
$\sigma$ fixes all the elements generated by $R$ because the operations 
are permutation invariant. Thus if $a\in A'$ and $s\in a\cap R$ then $R\subseteq a$
showing that $R$ is an atom of $\A'$. Since $\tau$ is a bijection, it follows that ${\sf s}_{\tau}R$ is also an atom of $\A$
and, it is easy to see that  all these atoms are pairwise disjoint.
That is if $\tau_1\neq \tau_2$, then ${\sf s}_{\tau_1}R\cap {\sf s}_{\tau_2}R=\emptyset$.
We now split each ${\sf s}_{\tau}R$ into abstract atoms ${\sf s}_{\tau}R_j$, $j\leq m$ and $\tau\in G_n$.
Let $(R_j:j\leq m)$ be a set of $m+1$ distinct elements, and let $\A_{k,n}$ be an algebra
such that
\begin{enumerate}
\item $\A'\subseteq \A_{k,n},$ the Boolean part of $\A_{k,n}$ is a Boolean algebra,
\item $R=\sum\{R_j:j\leq m\},$
\item ${\sf s}_{\tau}R_j$ are pairwise distinct atoms of $\A_k$ for each $\tau\in G_n$ and $j\leq m$ and 
${\sf c}_i{\sf s}_{\tau}R_j={\sf c}_i{\sf s}_{\tau}R$ for all $i<\omega$ and all $\tau\in G_n,$
\item each element of $\A_{k,n}$ is a join of element of $\A'$ and of some ${\sf s}_{\tau}R_j$'s,
\item ${\sf c}_i$ distributes over joins,
\item The ${\sf s}_{\tau}$'s are Boolean endomorphisms such that ${\sf s}_{\tau}{\sf s}_{\sigma}a={\sf s}_{\tau\circ \sigma}a$.
\end{enumerate}

The existence of such algebra is easy to show; furthermore they are unique up to isomorphim, see \cite{Andreka}, the comment right 
after the 
definition on p.168. Now we show that $\Rd_{ca}\A_{k,n}$ cannot be representable. This part of the proof is identical to 
Andr\'eka's proof but we include it for the sake of 
completeness. The idea is that we split $R$ into $m+1$ distinct atoms but $U_0$ has only $m$ elements, and those two 
conditions are incompatible in case there is a 
representation. The substitutions have to do with permuting the atoms and they do not contribute to this part of the proof.
For $i,j<\omega$, $i\neq j,$ recall that ${\sf s}_j^ix={\sf c}_i({\sf d}_{ij}\cdot x)$. Let 
$$\tau(x)=\prod_{i\leq m} {\sf s}_i^0{\sf c}_1\ldots {\sf c}_mx\cdot \prod_{i<j\leq m} 
-{\sf d}_{ij}$$
Then $\A'\models \tau(R)=0.$
Indeed we have  
$${\sf c}_1\ldots {\sf c}_mR={}^mU\times U_{m+1}\times \ldots $$
$${\sf s}_i^0{\sf c}_1\ldots {\sf c}_mR=U\times \ldots U_0\times 
U\times U_{m+1}\ldots $$
$$\bigcap {\sf s}_i^0{\sf c}_1\ldots {\sf c}_mR={}^{m+1}U_0\times U_{m+1}\times .$$
Then by $|U_0|\leq m$ there is no repitition free sequence in $^{m+1}U_0$. Thus
as claimed $\A'\models \tau(R)=0$.
Then $\A_{k,n}\models \tau(R)=0$. Assume that $\A_{k,n}$ is represented somehow.
Then there is a homomorphism 
$h:\A_{k,n}\to \langle \B(^{\omega}W), {\sf c}_i, {\sf d}_{ij}\rangle_{i,j<\omega}$ 
for some set $W$ such that $h(R)\neq \emptyset$.

By $h(R)\neq \emptyset$ there is some $s\in h(R)$. By $R\leq {\sf c}_0R_i$ we have
$h(R)\subseteq h(R_i)$, so there is a $w_i$ such that $s(0|w_i)\in h(R_i)$ 
for all $i\leq m$.
These $w_i$'s are distinct since the $R_i$'s are pairwise disjoint (they are distinct atoms)
and so are the 
$h(R_i)$'s.
Consider the sequence
$$z=\langle w_0, w_1, \ldots w_m, s_{m+1}, \ldots \rangle.$$
We show that $z\in \tau(h(R))$. Indeed let $i, j\leq m$, $i\neq j$, then $z\in -{\sf d}_{ij}$
by $w_i\neq w_j$. Next we show that $z\in {\sf s}_i^0{\sf c}_1\ldots {\sf c}_mh(R)$. 
By definition, $\langle w_i, s_1\ldots \rangle\in h(R_i)\subseteq h(R)$
so $\langle w_i, w_1,\ldots w_m, s_{m+1}, \rangle\in {\sf c}_1\ldots {\sf c}_mh(R)$ 
and thus $z\in {\sf c}_0({\sf d}_{0i}\cap {\sf c}_1\ldots {\sf c}_mh(R))={\sf s}_i^0{\sf c}_1\ldots {\sf c}_mh(R)$.
This contradicts that $\A_{k,n}\models \tau(R)=0$.

Next we show that the $k$ generated subalgebras of $\A_{k,n}$ are representable.
Let $G$ be given such that $|G|\leq k$. The idea is to use $G$ and define a ``small" subalgebra of $\A_{k,n}$ 
that contains $G$
and is representable. 
Define $R_i\equiv R_j$ iff
$$(\forall g\in G)(\forall \tau\in G_n)({\sf s}_{\tau}R_i\leq g\Longleftrightarrow {\sf s}_{\tau}R_j\leq g).$$
This is similar to the equivalence relation defined by Andr\'eka \cite{Andreka} p. 157; 
the difference is that substitutions have to come to the picture \cite{Andreka}p.189.
Then $\equiv$ is an equivalence relation on $\{R_j: j\leq m\}$ which has $\leq 2^{k.n!}$ blocks by $|G|\leq k$ and $G_n=n!.$
Let $p$ denote the number of blocks of $\equiv$, that is $p=|\{R_j/\equiv:j\leq m\}|\leq 2^{k.n!}\leq m$.
Now that $R$ is split into $p< m+1$ atoms, the incompatibility condition above no longer holds. 
Indeed, let
$$B=\{a\in A_{k,n}: (\forall i,j\leq m)(\forall \tau\in G_n)(R_i\equiv R_j\text { and }{\sf s}_{\tau}R_i\leq a\implies {\sf s}_{\tau}R_j\leq a\}.$$

We first show that $B$ is closed under the operations of $\A_{k,n}$, then we show 
that, unlike $\A_{k,n}$,  $B$ is the universe of a representable algebra. Let $i<l<\omega$
Clearly $B$ is closed under the Boolean operations. The diagonal element
${\sf d}_{il}\in \B$ since ${\sf s}_{\tau}R_j\nleq {\sf d}_{il}$ for all $j\leq m$ and $\tau\in G_n$.
Also $A'\subseteq B$ since ${\sf s}_{\tau}R$ is an atom of $\A'$ and ${\sf c}_ia\in A'$ for all $a\in A_{k,n}$.
Thus ${\sf c}_ib\in B$ for all $b\in B$.
Assume that $a\in B$ and let $\tau\in G_n$. Suppose that $R_i\equiv R_j$ and ${\sf s}_{\sigma}R_i\leq {\sf s}_{\tau}a$. Then 
${\sf s}_{\tau}{\sf s}_{\sigma}R_i\leq a$, so ${\sf s}_{\tau\circ \sigma}R_i\leq a$. Since $a\in B$ we get
that ${\sf s}_{\tau\circ \sigma}R_j={\sf s}_{\tau}{\sf s}_{\sigma}R_j\leq a$, and so ${\sf s}_{\sigma}R_j\leq{\sf s}_{\tau}a$.
Thus $B$ is also closed under substitutions.
Let $\B\subseteq \A_{k,n}$ be the subalgebra of $\A_{k,n}$ with universe $B$. Since $G\subseteq B$
it suffices to show that $\B$ is representable.
Let $\{y_j;j<p\}=\{\sum(R_j/\equiv):j\leq m\}$. 
Then $\{y_j:j<p\}$ is a partition of $R$ in $\B,$ ${\sf c}_iy_j={\sf c}_iR$ for all $j<p$ and 
$i<\omega$ and every element of $\B$ is a join of some element of $\A'$ and of finitely many of 
${\sf s}_{\tau}y_j$'s. Recall that $p\leq m$.
We now split $R$ into $m$ `real' atoms, cf. \cite{Andreka} p.167, lemma 2. We define an equivalence relation on $R$. For any $s,z\in R$
$$s\sim z\Longleftrightarrow |\{i\in \omega: s_i\neq z_i\}|<\omega.$$
Let $S\subseteq R$ be a set of reprsentatives of $\sim$. 
Consider the group $Z_m$ of integers modulo $m$. 
(Any finite abelian group with $m$ elements 
will do.) 
For any $s\in S$ and $i\in \omega$ let $f_i^s:U_i\to Z_m$ be an onto map such that $f_i^s(s_i)=0$.
For $j<m$ define
$$R_j^s=\{z\in R: \sum\{f_i^s:i\in \omega\}=j\}$$
and 
$$R_j"=\bigcup \{R_j^s: s\in S\}.$$
Then $\{R_0'', \ldots  R_{m-1}''\}$ is a partition of $R$ such that ${\sf c}_iR_j''= {\sf c}_iR$ for all $i<\omega$ and $j<m$.
Let $\A''$ be the subalgebra of 
$\langle \B(^{\omega}U), {\sf c}_i, {\sf d}_{ij}, {\sf s}_{\tau}\rangle_{i,j<\omega, \tau\in G_n}$
generated by $R_0'',\ldots R_{m-1}''.$
Let 
$$\R=\{{\sf s}_{\sigma}R_j'':\sigma\in G_n, j<m\}.$$
Let $$H=\{a+\sum X: a\in A', X\subseteq_{\omega}\R\}.$$
Clearly $H\subseteq A''$ and $H$ is closed under the boolean operations. Also
because transformations considered are bijections we have
$${\sf c}_i{\sf s}_{\sigma}R_j={\sf c}_i{\sf s}_{\sigma}R\text { for all $j<m$ and }\sigma\in G_n.$$
Thus $H$ is closed under ${\sf c}_i.$ Also $H$ is closed under substitutions.
Finally ${\sf d}_{ij}\in A'\subseteq H.$
We have proved that $H=A''$. This implies that every element of $\R$ is an atom of $\A''$.
We now show that $\B$ is embeddable in $\A''$, and hence will be representable.
Define for all $j<p-1$,
$$R_j'=R_j'',$$
and
$$R_{p-1}'=\bigcup\{R_j'': p-1\leq j<m\}$$
Then define for $b\in B$:
$$h(b)=(b-\sum_{\tau\in G_n}{\sf s}_{\tau}R)\cup\bigcup\{{\sf s}_{\tau}R_j': \tau\in G_n,  j<p, {\sf s}_{\tau}y_j\leq b\}.$$
It is clear that $h$ is one one, 
preseves the Boolean operations and the diagonal elements and is the identity on $A'$.
Now we check cylindrifications and substitutions.
\begin{equation*}
\begin{split}
{\sf c}_ih(b)&={\sf c}_i[(b-\sum {\sf s}_{\tau}R)\cup\bigcup\{{\sf s}_{\tau}R_j': \tau\in G_n, j<p, {\sf s}_{\tau}y_j\leq b\}]\\
&={\sf c}_i(b-\sum_{\tau\in G_n} {\sf s}_{\tau}R)\cup\bigcup\{{\sf c}_i{\sf s}_{\tau}R_j':\tau\in G_n, j<p, {\sf s}_{\tau}y_j\leq b\}\\
&={\sf c}_i(b-\sum_{\tau\in G_n} {\sf s}_{\tau}R\cup\bigcup\{{\sf c}_i{\sf s}_{\tau}y_j: \tau\in G_n, j<p, {\sf s}_{\tau}y_j\leq b\}\\
&={\sf c}_i[(b-\sum_{\tau\in G_n} {\sf s}_{\tau}R)\cup\bigcup\{{\sf s}_{\tau}y_j,\tau\in G_n, j<p, {\sf s}_{\tau}y_j\leq b\}]\\
&={\sf c}_ib\\
\end{split}
\end{equation*}
On the other hand
$$h{\sf c}_i(b)=({\sf c}_ib-\sum_{\tau\in G_n} {\sf s}_{\tau}R)\cup\bigcup\{{\sf s}_{\tau}R_j': {\sf s}_{\tau}y_j\leq {\sf c}_ib\}={\sf c}_ib.$$
Preservation of substitutions follows from the fact that the substitutions are Boolean endomorphisms. In more detail, let $\sigma\in G_n$, then:
\begin{equation*}
\begin{split}
{\sf s}_{\sigma}h(b)&={\sf s}_{\sigma}[(b-\sum_{\tau\in G_n} {\sf s}_{\tau}R)\cup\bigcup\{{\sf s}_{\tau}R_j':\tau\in G_n, j<p, {\sf s}_{\tau}y_j\leq b\}]\\
&=({\sf s}_{\sigma}b-\sum_{\tau\in G_n} {\sf s}_{\sigma}{\sf s}_{\tau}R)\cup\bigcup\{{\sf s}_{\sigma}
{\sf s}_{\tau}R_j': \tau\in G_n, j<p, {\sf s}_{\tau}y_j\leq b\}]\\
&=({\sf s}_{\sigma}b-\sum_{\tau\in G_n} {\sf s}_{\sigma\circ \tau}R)\cup\bigcup\{{\sf s}_{\sigma\circ \tau}R_j': \tau\in G_n, j<p, {\sf s}_{\tau}y_j\leq b\}]\\
&=({\sf s}_{\sigma}b-\sum_{\tau\in G_n} {\sf s}_{\tau}R)\cup\bigcup\{{\sf s}_{\tau}R_j':  \tau\in G_n, j<p, {\sf s}_{\tau}y_j\leq b\}]\\
\end{split}
\end{equation*}

But for fixed $\sigma,$ we have $\{\sigma\circ \tau: \tau\in G_n\}=G_n$ and so
$${\sf s}_{\sigma}h(b)=h({\sf s}_{\sigma}(b)).$$

For every $k,n<\omega$ we have constructed an algebra $\A_{k,n}$ such that 
$\Rd_{\CA}\A_{k,n} \notin \RCA_{\omega}$ and the $k$-generated subalgebras of 
$\A_{k,n}$ are representable.  
We should point out that the  ``finite dimensional version" of the $\A_{k,n}$'s 
were constructed in \cite{c}, and their construction can be recovered from 
the proof of Theorem 6 in \cite{Andreka} which addresses the finite dimensional case but in a more general 
setting allowing arbitrary unary additive permutation invariant operations expanding 
those of $\RCA_n$. We note that the latter result does not survive the infinite dimensional case. There are easy examples, cf. \cite{Andreka} p.192 
and \cite{ex}. 

\item We show that $\A_{k,n}$ has a representation which preserves all operations except for the diagonal elements. 
That is, its quasipolyadic reduct is 
representable.
The proof is analogous to that of Andr\'eka's on of Claim 16 on p.194 of \cite{Andreka}.
Let $\A_{k,n}$ be the algebra obtained by splitting the atom $R$ in $\A'$ as in the above proof. Then
$\A_{k,n}$ is not representable,  but its $k$ generated subalgebras are representable.
We show that there is a representation of $\A_{k,n}$ in which all operations are preserved except for the diagonal elements.
Let $U_i$, $i<\omega$ be a sequence of pairwise disjoint sets such that $|U_0|=m\geq 2^{k.n!+1}$  and $|U_i|\geq m+1$. Let $R$ be as above
except that it is defined via the new $U_i$'s.   Let $(R_j:j\leq m)$ be the splitting of $R$ in $\A_{k,n}$. 
Let $W\supset U$ (properly).
Let $W_0=U_0\cup (W\sim U)$, and $W_i=U_i$ for $0<i<\omega$.
First we define a function $h:\wp(^{\omega}U)\to \wp(^{\omega}W)$ with the desired properties and $h(R)=\prod_{i<\omega}W_i$.
Let $t: W\to U$ be a surjective function which is the identity on $U$ and which maps $W_0$ to $U_0$. Define $g:{}^{\omega}W\to {}^{\omega}U$ 
by $g(s)=t\circ s$ for all $s\in {}^{\omega}W$ and for all $x\subseteq {}^{\omega}U,$ define
$$h(x)=\{s\in {}^{\omega}U: g(s)\in x\}.$$
Since $|W_i|\geq m+1$ for all $i<\omega$, the incompatibity condition between the number of atoms 
splitting $R$ and the number of elements in $|W_0|$
used in the representation vanishes, 
so there is a real partition $(S_j: j\leq m)$ 
of $S=\prod_{i<\omega}W_i$ such that ${\sf c}_iS_j={\sf c}_iS$ for all $ i<\omega$ and $j\leq m$.
Then $({\sf s}_{\sigma}S_j: j\leq m)$ is an analogous partition of ${\sf s}_{\sigma}S$ for $\sigma\in G_n$.
Let $X_{\sigma,j}={\sf s}_{\sigma}^{\A_{k,n}}R_j$ for $j\leq m$.
Define $\bar{h}:\A_{k,n}\to \wp(^{\omega}W)$ by
$$\bar{h}(a)=h(a), \ \ a\in A'$$
$$\bar{h}(X_{\sigma j})={\sf s}_{\sigma}S_j, \sigma \in G_n, j\leq m$$
and
$$\bar{h}(x+y)=\bar{h}(x)+\bar{h}(y), x,y\in A_{k,n}.$$
It is easy to check using lemma  (iv) in \cite{Andreka} that $\bar{h}$ is as desired. In fact $\bar{h}$ 
preserves all the quasipolyadic operations including substitutions 
corresponding to replacements, which are now no longer definable, because we have discarded diagonal elements.
The reasoning is as follows \cite{Andreka} p.194. 
For $i,j\in n$, 
the quantifier free formula $x\leq -{\sf d}_{ij}\to {\sf s}_i^jx$ is valid in representable algebras 
hence it is valid in $\A_{k,n}$ since its  $k$ generated subalgebra are representable. 
Let $\sigma\in G_n$, $l\leq m$. Then ${\sf s}_i^j(X_{\sigma, l})=0$ in $\A_{k,n}$.
Now
$$\bar{h}({\sf s}_i^j(X_{\sigma, l}))=\bar{h}(0)=0={\sf s}_i^j {\sf s}_{\sigma}S_j={\sf s}_i^jh(X_{\sigma l}).$$ 
Assume that $a\in A_{k,n}$. Then
$$\bar{h}({\sf s}_i^ja)=h({\sf s}_i^ja)={\sf s}_i^jh(a)={\sf s}_i^j\bar{h}(a).$$ 
Since both $\bar{h}$ and ${\sf s}_i^j$ are additive we get the required. 
\end{enumarab}

{\bf Part II}
\begin{enumarab}
\item Here is where we really start the non-trivial modification of Andr\'eka's splitting. 
For $n\in \omega$ and $m=2^{k.n!+1}$, we denote $\A_{k,n}$ by $split(\A', R, m, n)$.
This is perfectly legitimate since the algebra $\A_{k,n}$ is determined uniquely by $R,$ $\A'$, $m$ and $n$. 
Recall that $m$ is the number of atoms splitting $R$, while $n$ is the finite number of substitutions available.
For $n_1<n_2$, we denote by $\Rd_nsplit(\A', R, m, n_2)$ 
the reduct of $split(\A', R, m, n_2)$ obtained by restricting substitutions to $G_{n_1}.$ 
Let $m_1<m_2$ and $n_1<n_2$. Then we claim that
$$split(\A', R, m_1,n_1)\text { embeds into }\Rd_{n_1}split(\A', R, m_2, n_2).$$
This part of the proof is analogous to Andr\'eka's proofs in \cite{Andreka}, lemma 3, on splitting elements in cylindric algebras.
Indeed, let $$\chi: m_1\to m_2$$ 
be such that the set $\chi(j),$ $j<m_1$ are non empty and pairwise disjoint, and
$$\bigcup \{\chi(j):j<m_1\}=m_2.$$
For $x\in split(\A, R, m_1, n_1),$ let 
$$J_{\tau}(x)=\{j<m_1: {\sf s}_{\tau}R_j\leq x\}.$$
Let $(R_i: i\leq m_2)$ be the splitting of $R$ in $split(\A', R, m_2, n_2)$.
Define 
$$h(x)=(x-\sum {\sf s}_{\tau}R)+\sum\{{\sf s}_{\tau}R_i: \tau\in G_{n_1}, i\in \bigcup \{\chi(j): j\in J_{\tau}(x)\}\}.$$
Here we are considering $G_{n_1}$ as a subset of $G_{n_2}$.
It is easy to check that $h(x)$ is a Boolean homomorphism and that $h(x)\neq 0$ whenever $0\neq x\leq {\sf s}_{\tau}R,$ for $\tau\in G_{n_1}$.
Thus $h$ is one to one.
Let $i\in \omega$ and $x\in split(\A, R, m_1, n_1)$. If
$x\cdot {\sf s}_{\tau}R=0$ for all $\tau,$ then $x\in A'$, hence $h({\sf c}_ix)={\sf c}_ih(x)$. 
So assume that there is a $\tau\in G_n$ such that $x\cdot {\sf s}_{\tau}R\neq 0$. Then ${\sf c}_i(x\cdot {\sf s}_{\tau}R)={\sf c}_i{\sf s}_{\tau}R$ and 
${\sf c}_ih(x\cdot {\sf s}_{\tau}R))={\sf c}_i{{\sf s}_{\tau} R}$ by $0\neq h(x\cdot {\sf s}_{\tau}R)\leq {\sf s}_{\tau}R$.
Now
\begin{equation*}
\begin{split}
h({\sf c}_ix)&=h({\sf c}_i(x-{\sf s}_{\tau}R)+{\sf c}_i(x\cdot {\sf s}_{\tau}R))\\
&={\sf c}_i(x-{\sf s}_{\tau}R)+{\sf c}_iR.\\
\end{split}
\end{equation*}
\begin{equation*}
\begin{split}
{\sf c}_ih(x)&={\sf c}_i(h(x-{\sf s}_{\tau}R)+x\cdot {\sf s}_{\tau}R))\\
&={\sf c}_i(h(x-{\sf s}_{\tau}R))+h(x\cdot {\sf s}_{\tau}R)\\
&={\sf c}_i(x-{\sf s}_{\tau}R)+{\sf c}_iR.\\
\end{split}
\end{equation*}
We have proved that 
$${\sf c}_ih(x)=h({\sf c}_ix).$$ Now we turn to substitutions. Let $\sigma\in G_{n_1}$. Then we have
\begin{equation*}
\begin{split}
{\sf s}_{\sigma}h(x)&={\sf s}_{\sigma}[(x-\sum_{\tau\in G_{n_1}} {\sf s}_{\tau}R)+\sum\{{\sf s}_{\tau}R_i, \tau\in G_{n_1}, i\in \bigcup\{\chi(j): j\in J_{\tau}x\}\})]\\
&=({\sf s}_{\sigma}x-\sum_{\tau\in G_{n_1}}{\sf s}_{\sigma\circ \tau}R)+
\sum\{{\sf s}_{\sigma\circ \tau}R_i:\tau\in G_{n_1}, i\in \bigcup\{\chi(j): j\in J_{\tau}x\}\}.\\
\end{split}
\end{equation*}
Since $\{\sigma\circ \tau:\tau\in G_{n_1}\}=G_{n_1}$ then we have:
$${\sf s}_{\sigma}h(b)=h({\sf s}_{\sigma}(b).)$$

\item We have a sequence of algebras $(\A_{k,i}:i\in \omega\sim{0})$ such that for $n<m,$ 
we can assume by the embeddings proved to exist in the previous item
that $\A_{k,n}$ is a subreduct (subalgebra of a reduct) 
of $\A_{k,m}$.  
Form the natural direct limit of such algebras which is the (reduct directed) union
call it $\A_k$. That is $A_k=\bigcup_{n\in \omega}A_{k,n}$, and the operations are defined the obvious way.
For example if $i<\omega$, and $a\in A_k$, then $i\in n$ and $a\in A_{n,k}$ for some $n$; set ${\sf c}_i^{\A_k}a={\sf c}_i^{\A_{k,n}}a$.
These are well defined. The other operations are defined analogously, where we only define the ${\sf s}_{[i,j]}$'s for $i,j\in \omega$. 
Clearly, $\Rd_{ca}\A_k$ is not representable, for else $\Rd_{ca}\A_{k,n}$ would be representable for all $n\in \omega$.

\item Let $|G|\leq k$. Then $G\subseteq \A_{k,n}$ for some $n$. If $\Sg^{\A_{k}}G$ is not representable then there exists 
$l\geq n$ such that $G\subseteq \A_{k,l}$ and 
$\Sg^{\A_{k,l}}G$ is not representable, contradiction.
To see this we can show directly that $\Sg^{\A_k}G$ has to be representable. 
We show that every equation $\tau=\sigma$  valid in the variety
$\RQEA_{\omega}$ is valid in $\Sg^{\A_k}G$. 
Let $v_1,\ldots v_k$ be the variables occuring in this equation, and let $b_1,\ldots b_k$ be arbitrary elements of $\Sg^{\A_k}G$. 
We show that $\tau(b_1,\ldots b_k)=\sigma(b_1\ldots b_k)$. Now there are terms 
$\eta_1\ldots \eta_k$ written up from elements of $G$ such that $b_1=\eta_1\ldots b_k=\eta_k$, then we need to show that 
$\tau(\eta_1,\ldots \eta_k)=\sigma(\eta_1, \ldots \eta_k).$ 
This is another equation written up from elements of $G$, which is also valid in $\RQEA_{\omega}$. 
Let $n$ be an upper bound for the indices occuring in this equation and let $l>n$ be such that $G\subseteq \A_{k,l}$. 
Then the above equation is valid in $\Sg^{\Rd_n\A}G$ since the latter is representable. 
Hence the equation $\tau=\sigma$ holds in $\Sg^{\A_k}G$ at the evaluation $b_1,\ldots b_k$ of variables.

\item Let $\Sigma_n^v$ be the set of universal formulas using only $n$ substitutions and $k$ variables valid in $\RQEA_{\omega}$, 
and let $\Sigma_n^d$ be the set of universal formulas
using only $n$ substitutions and no diagonal elements valid in $\RQEA_{\omega}$.  By $n$ substitutions we understand the set 
$\{{\sf s}_{[i,j]}: i,j\in n\}.$
Then $\A_{k,n}\models \Sigma_n^v\cup \Sigma_n^d$. $\A_{k,n}\models \Sigma_n^v$ because the $k$ generated subalgebras
of $\A_{k,n}$ are representable, while $\A_{k,n}\models \Sigma_n^d$ because $\A_{k,n}$ has a representation that preserves all operations except
for diagonal elements.  Indeed, let $\phi\in \Sigma_n^d$, then there is a representation of $\A_{k,n}$ in which all operations 
are the natural ones except for the diagonal elements. 
This means that (after discarding the diagonal elements) there is a one to one homomorphism 
$h:\A^d\to \P^d$ where $\A^d=(A_{k,n}, +, \cdot , {\sf c}_k, {\sf s}_{[i,j]}, {\sf s} _i^j)_{k\in \omega, i,j\in n}\text { and } 
\P^d=(\B(^{\omega}W), {\sf c}_k^W, {\sf s}_{[i,j]}^W, {\sf s}_{[i|j]}^W)_{k\in \omega, i,j\in n},$ 
for some infinite set $W$. 
Now let $\P=(\B(^{\omega}W), {\sf c}_k^W, {\sf s}_{[i,j]}^W, {\sf s}_{[i|j]}^W, {\sf d}_{kl}^W)_{k,l\in \omega, i,j\in n}.$
Then we have that $\P\models \phi$ because $\phi$ is valid 
and so  $\P^d\models \phi$ due to the fact that  no diagonal elements  occur in $\phi$. 
Then $\A^d\models \phi$ because $\A^d$ is isomorphic to a subalgebra of $\P^d$ and $\phi$ is quantifier free. Therefore 
$\A_{k,n}\models \phi$.
Let $$\Sigma^v=\bigcup_{n\in \omega}\Sigma_n^v 
\text { and }\Sigma^d=\bigcup_{n\in \omega}\Sigma_n^d$$ 
Hence $\A_k\models \Sigma^v\cup \Sigma^d.$ For if not then there exists a quantifier free  formula $\phi(x_1,\ldots x_m)\in \Sigma^v\cup \Sigma^d$,
and $b_1,\ldots b_m$ such that $\phi[b_1,\ldots b_n]$ does not hold in $\A_k$. We have $b_1\ldots b_m\in \A_{k,i}$ for some $i\in \omega$. 
Take $n$ large enough $\geq i$ so that
$\phi\in \Sigma_n^v\cup \Sigma_n^d$.   Then $\A_{k,n}$ does not model $\phi$, a contradiction.
Now let $\Sigma$ be  a set of quantifier free formulas axiomatizing  $\RQEA_{\omega}$, then $\A_k$ does not model $\Sigma$ since $\A_k$ is not 
representable, so there exists a formula $\phi\in \Sigma$ such that
$\phi\notin \Sigma^v\cup \Sigma^d.$ Then $\phi$ contains more than $k$ variables and a diagonal constant occurs in $\phi$.
\end{enumarab}
\end{demo}

We immediately get the following answer to Andr\'eka's question formulated on p. 193 of \cite{Andreka}.

\begin{athm}{Corollary 2} The variety $\RQEA_{\omega}$ is not axiomatizable over $\RQA_{\omega}$ 
with a set of universal formulas containing infinitely many variables.
\end{athm}
One can show, using the modified method of splitting here, that all theorems in \cite{Andreka} on complexity of axiomatizations 
generalize to $\RQEA_{\omega}$ with the sole exception of Theorem 5, 
which is false for infinite dimensions \cite{ex}. As a sample we give the following theorem 
which can be proved by some modifications of the cited theorems
in the proof; this modifications are not hard, most of them can be found in \cite{c}, proof of theorem 3.1. The basic idea 
in the proof is to show that certain cylindric homomorphisms
(between cylindric algebras) remain to be quasipolyadic equality algebra homomorphisms 
(between their corresponding natural quasipolyadic equality expansions), that is when we add substitutions corresponding to transpositions.

\begin{athm}{Theorem 3} Let $\Sigma$ be a set of equations 
axiomatizing $\RQEA_{\omega}.$  Let $l,k, k' <\omega$. 
Then $\Sigma$ contains
infinitely equations in which $-$ occurs, one of $+$ or $\cdot$ occurs  a diagonal or a permutation with index $l$ occurs, more
than $k'$ cylindrifications and more than $k$ variables occur.
\end{athm}
\begin{demo}{Sketch of Proof} Let $n,k\in \omega\sim\{0\}$. Let $\A_{k,n}$ be the non-representable algebra constructed above 
obtained by splitting $\A'$ into $m\geq 2^{k.n!+1}$ atoms and we require that $m\geq k'$ as well. One then shows that 
the complementation free reduct  $\A_{k,n}^{-}$ of $\A_{k,n}$ is a homomorphic image of a subalgebra $\C$ of 
the complemention free reduct of $\P^{-}$ of a  a representable $\P$, \cite{Andreka}, cf. Theorem 7, p.163. The algebra
$\A_{k,n}$ can be represented such that every operation except for $\cup$ and $\cap$ are the natural ones, cf. \cite{Andreka} p.200 and \cite{c} for the necessary 
modifications.
For any $I\subseteq \omega$, $|I|=m$ there is an infinite set $W$ an an embedding  
from $\A_{k,n}\to (\B(^{\omega}W), {\sf c}_i, {\sf d}_{ij}, {\sf s}_{\tau})_{i,j\in \omega, \tau\in G_n}$ 
which is a homomorphism with respect to all operations of $\A_{k,n}$ except for ${\sf c}_i$ $i\notin I$, cf. \cite{Andreka} p.172, Theorem 3.
One just has to show that the map $h:\A'\to \wp(^{\omega}W)$ defined on p. 174 preserves substitutions, which is 
straightforward from the definition of the map $g$ defined on 
p.173.
There is an infinite set $W,$ such that there is an embedding $h:\A\to (\B(^{\omega}W), {\sf c}_i, {\sf d}_{ij}, {\sf s}_{\tau})_{i,j<n, \tau\in G_n}$ 
such that $h$ is a homomorphism preserving all operations except for ${\sf d}_{il}$ and ${\sf s}_{[i,l]}$ if $i,l\in n$, cf. p.176 Claim 6. 
This will prove the theorem because of the following reasoning.  
By $n$ substitutions we understand the set $\{{\sf s}_{[i,j]}: i,j\in n\}$. 
Let $\Sigma_n^-$ denote the set of equations  without complementation in which only $n$ substitutions occur, 
$\Sigma_n^v$ be the set of equations
which contains at most $k$ variables in which only $n$ substitutions occur, $\Sigma_n^c$ be the set of equations 
in which only $k'$ cylindrifications and $n$ substitutions occur,
$\Sigma_n^{ds}$ be the set of equations in which at most $n$ substitutions occur and no diagonal nor substitutions with index $l$ occurs, 
and $\Sigma_n^{Bool}$ the set of equations 
that does not contain $\cdot$ nor  $+$ and $n$ substitutions occur, all valid in $\RQEA_{\omega}$.  Then 
$\A_{k,n}\models \Sigma_n^-\cup \Sigma_n^v\cup \Sigma_n^c\cup \Sigma_n^{ds}\cup \Sigma_n^{Bool}.$ 
Indeed, the algebra $\A_{k,n}\models \Sigma_n^-$ because of the following reasoning.
Let $\C$ and $\P$ be as above. Then $\P\models \Sigma_n^-$ because $\P$ is representable. So $\P^-\models \Sigma_n^-$ 
because $-$ does not occur in 
$\Sigma_n^-$. Now  $\C\models \Sigma_n^-$ by 
$\C\subseteq \P^{-}$, and so $\A_{k,n}^-\models \Sigma_n^-$ since $\A^{-}_{k,n}$ is a homomorphic image of $\C$ and $\Sigma_n^-$ 
consists of equations. 
Then $\A_{k,n}\models \Sigma_n^-$. 
This together with previous reasoning proves that $\A_{k,n}\models  \Sigma_n$ where $\Sigma_n$ is the above (union) of formulas. 
In more detail $\A_{k,n}\models \Sigma^v_n$ because its $k$ generated subalgebras are representable, 
$\A_{k,n}\models \Sigma_n^{dp}$ because it has a representation that preserves all elements except for 
diagonals and substitutions with index $l$, and so forth. Then we can infer that  
$\A_k\models \bigcup_{n\in \omega} \Sigma_n=\Gamma.$
For if not, then we can choose $n$ large enough such that $\A_{k,n}$ does not model 
$\Sigma_n$. But $\A_k$ is not representable hence any equational 
axiomatization of the representable algebras contain a formula that is outside $\Gamma$.
Thus the required follows. 
\end{demo}
For axiomatization with universal formulas, using the same ideas as above, except for those involving complementation,
we obtain the slightly weaker:
\begin{athm}{Theorem 4}
Let $\Sigma$ be a set of quantifier free formulas 
axiomatizing $\RQEA_{\omega}.$  Let $l,k, k' <\omega$. 
Then $\Sigma$ contains
infinitely equations in which  one of $+$ or $\cdot$ occurs  a diagonal or a permutation with index $l$ occurs, more
than $k'$ cylindrifications and more than $k$ variables occur
\end{athm}
Our corollary 2 is  evidence that $\RQA_{\omega}$ can be axiomatized by an infinite set of universal formulas containing only 
finitely many variables. The investigations carried out in this paper are also strongly 
related to the so-called finitizability (very central) problem in algebraic logic, 
\cite{Bulletin}, \cite{SG}, \cite{Sain2000}, \cite{Fer} which roughly asks for nice algebraisations of (algebraisable) 
extensions of first order logic that avoid 
complex non-finite axiomatizability results that dominate 
standard algebraisations like representable cylindric algebras and quasipolyadic 
algebras (with and without equality). The paper deals with sharpening negative non-finite axiomatizabilty results.

However, rather surprisingly, positive results circumventing or side-stepping 
such non-finite axiomatizations with a high degree of complexity  for this problem do exist \cite{SG}, \cite{Sain2000}. 
But the presence of diagonal elements keeps the problem on the tough side and this suggests 
that the class of representable quasipolyadic algebras are easior to tame than that of quasipolyadic 
equality algebras by passing to finitary expansions via finitely many infinitary substitutions 
\cite {Sain2000}, \cite{SG}. This is compatible with the fact that the `distance' between 
the two is infinite, which is one reading of 
our main result. 

Another trend for obtaining positive solutions for the finitizability problem is 
to relativize representations \cite{Fer} obtaining Henkin-like semantics, with the slogan being relativization turns negative results to positive ones.
This trend originates with Henkin and Resek in the latter's Ph.D dissertation which was never published, cf. \cite{Fer} 
for reviewing the development of this very 
rich subject. The proof in \cite{Fer} also illustrates that generalizing results from the cylindric paradigm to the quasipolyadic equality one,
is far from trivial. Such an approach has close affinity with modal logic.

We note that a recursive axiomatization of $\RQEA_{\alpha},$ $\alpha$ an infinite ordinal is given in \cite{c} 
by using the methods of Hirsch and Hodkinson of synthesising 
axioms by games.

\end{document}